\magnification1200

\input eplain

\input epsf

\newcount\secno
\newcount\thmno
\def\here{\the\secno.\the\thmno}

\newcount\figno

\def\sec#1{\par\vskip8mm
\advance \secno by 1
\centerline{{\bf \the\secno\ #1}}
\thmno=0
\figno=0
\par\bigskip\noindent}

\def\thm#1#2{\par\bigskip\noindent
\advance\thmno by 1
{\bf #1 \here.} {\it #2}
\smallskip}

\def\thmnamed#1#2#3{\par\medskip\noindent
\advance\thmno by 1
{\bf #1 \here\ (#2).} {\it #3}
\smallskip}

\def\rmk#1{\par\medskip\noindent
\advance\thmno by 1
{\bf #1 \here. }}

\def\proof{\noindent{\it Proof: }}

\def\figure#1{\global\advance\figno by 1
\bigskip
\centerline{Figure~\the\secno.\the\figno: #1}}

\def\eps{\varepsilon}
\def\pr{{\rm pr}}

\secno=0

\everyfootnote={\leftskip=1cm \parindent=0cm 
                \rm}


\def\id{{\rm Id}}
\def\r{{\bf R}}
\def\z{{\bf Z}}

\def\cstr{contact structure}
\def\ot{overtwisted}
\def\ocstr{overtwisted \cstr}
\def\difmm{diffeomorphism}
\def\cmm{contactomorphism}

\def\dif{{\it Diff}}
\def\tdif{\widetilde{\it Diff}}
\def\cont{{\it Cont}}
\def\heq{homotopy equivalen}
\def\nbd{neighborhood}
\def\qed{{\hfill \boxit{}}}

\def\N{{\cal N}}

\def\rel{\mathop{\rm rel}\nolimits}


\def\dy{\cite{dymara:masters}}


\font\scaps=cmcsc10

\centerline{\bf THE GROUP OF CONTACTOMORPHISMS OF THE SPHERE}
\centerline{\bf FIXING AN OVERTWISTED DISK} 
\bigskip
\centerline{Katarzyna Dymara\footnote{*}{Partially supported by KBN grant
2 P03A 017 25.}}
\centerline{Institute of Mathematics}
\centerline{Wroc\l aw University}
\bigskip
\centerline{April 2005}
\vskip1cm

{\narrower
\centerline{\scaps Abstract}

We calculate the weak homotopy type of the group of contactomorphisms of
the three-sphere which coincide with the identity on (a neighborhood
of) an overtwisted disk.

}

\sec{Preliminaries. Known results}
\definexref{sec-prelim}{\the\secno}{section}

\noindent
A {\it \cstr} on a 3--dimensional manifold is a field of planes 
defined (at least locally) as the kernel of a 1--form 
$\alpha$ such that $\alpha\wedge d\alpha$ nowhere vanishes.

\smallskip
We say that two \cstr s $\zeta_0$ and $\zeta_1$ are 
{\it isotopic} if there exist a smooth family 
of \cstr s $\{\zeta_t, t\in[0,1]\}$.

We say that $\zeta_0$ and $\zeta_1$ are {\it contactomorphic}\/ if 
there exists a \difmm\ $f:M\to M$ such that 
$f_*(\zeta_0)=\zeta_1$; equivalently, if 
$\zeta_0=\ker\alpha_0$ and $\zeta_1=\ker\alpha_1$, then 
$(f^{-1})^*(\alpha_1)=t\alpha_0$ for some non-zero function $t:M\to\r$. 
Such $f$ is referred to as a {\it \cmm}.

Gray proved in \cite{gray} that on a closed manifold 
two contact structures are contactomorphic if and only if
they are isotopic.
His proof (which uses a vector field whose flow consists of 
the desired  \difmm s, constructed locally 
and glued together by means of a partition of unity) 
can be applied---without essential changes---to a relative situation, 
either in the sense of considering \cstr s on a manifold modulo
a compact set, or fixing the \cstr\ along a (contractible)
subset of the parameter space. Here we will simply state
Gray's theorem in necessary generality.


\thmnamed{Theorem}{Gray's Theorem}{Let $\{\zeta_t, t\in D^n\}$ be 
a smooth family of contact structures on a closed manifold $M$.
Assume that $\zeta_t|_A=\zeta_{t_0}|_A$ for a compact set $A\subset M$ 
and for all $t\in D^n$. Moreover, let $\zeta_t=\zeta_{t_0}$ for  all  
$t\in D'$, where $D'$ is a contractible subset of $D^n$. 
Then there exists a family $\{\phi_t, t\in D^n\}$ of \difmm s 
$\phi_t:M\to M$ such that for all $t$, $\phi_*\zeta_t=\zeta_{t_0}$, 
for all $t\in D^n$ $\phi_t|_A=\id_A$ and for all $t\in D'$ $\phi_t=\id_M$.}
\definexref{thm-gray2}{\here}{thm}


\medskip
A \cstr\ is called \ot\ if it contains a 2-dimensional disk which is tangent
to the \cstr\ along boundary. For example, consider the \cstr\ $\xi$ on
$\r^3$ defined as the kernel of the 1-form (written in cylindrical coordinates
$r,\theta,z$) 
$$
\cos r\,dz-\varrho(r)\,\sin r\,d\theta,
$$ 
where $\varrho:\r_+\to\r_+$ is a smooth function such that $\varrho(0)=0$, 
$\varrho'(0)>0$, $\varrho'(r)\geq0$ for all $r\in\r$
(we introduce the function $\varrho$ simply to make sure that 
the form is well-defined and smooth at $r=0$).

The disk $\Delta=\{(r,\theta,z) : r\le\pi, z=0\}$ is indeed
tangent to the \cstr\ along boundary.
Actually, every \ocstr\ contains a contactomorphic copy of (a \nbd\ of) $\Delta$;
this copy is referred to as an {\it \ot\ disk}.

A \cstr\ which is not \ot\ is called {\it tight}.


\thmnamed{Theorem}{Eliashberg's Theorem}{For $\Delta$ a 2-dimensional disk 
in an arbitrary 3-manifold $M$,
let $\xi_\Delta$ be a \cstr\ on a \nbd\ of $\Delta$ 
for which $\Delta$ is an \ot\ disk.
Denote by $\cont(M\rel\Delta)$ the space of \cstr s on $M$
which coincide with $\xi_\Delta$ on a \nbd\ of $\Delta$.
Moreover, let ${\it Distr}(M\rel\Delta)$ be the space of all plane
distributions on $M$ which coincide with $\xi_\Delta$ on a \nbd\ of $\Delta$.
Then the natural embedding
$\cont(M\rel\Delta)\to{\it Distr}(M\rel\Delta)$ is a weak \heq ce. }
\definexref{thm-eliashberg}{\here}{thm}


\smallskip
On $\pi_0$ level this means that isotopy classes of 
\cstr s \ot\ along a fixed disk (which, since embeddings of a disk
into a connected manifold are all isotopic, actually exhaust 
all isotopy classes of \ocstr s) remain in a one-to-one correspondence
with homotopy classes of plane fields.

to understand better this space, 
let us fix a parallelization of the manifold $M$, 
i.e.~a triple of vector fields $(u,v,w)$ such that 
$\bigl(u(x),v(x)w(x)\bigr)$ forms a basis of $T_xM$. 
Using this parallelization, we can identify a co-oriented 
plane field $\zeta$ with its normal Gauss map 
$G_\zeta:M\to S^2$ as follows: 
let $n_\zeta(x)$ be the unit normal vector to $\zeta$ at $x$, then 
$G_\zeta(x)=(n_1,n_2,n_3)$ for 
$n_\zeta(x)=n_1u(x)+n_2v(x)+n_3w(x)$.
For $S^3$, which is our main object of interest,
the space of homotopy classes of co-oriented\numberedfootnote{Of course, 
all plane fields on $S^3$
(or any manifold with trivial $H_1$) are co-orientable; as a matter of convenience,
we regard them as pre-equipped with one of the two possible co-orientations.}
plane fields (and hence also the space of isotopy classes of \ocstr s) 
is parameterized by 
$\pi_0\bigl({\it Map}(S^3{\to}S^2)\bigr)=\pi_3S^2\simeq{\bf Z}$.

Moreover, we may consider a parallelized plane field $\zeta$,
i.e.~one equipped with its own trivialization as a bundle.
This yields another Gauss map, $\widetilde G_\xi:M\to SO(3)$,
cf.~\dy, p.~301.
Parallelization and co-orientation of $\zeta$ form together 
a trivialization of the whole tangent bundle for which both Gauss maps
$G_\zeta$ and $\widetilde G$ are constant.

\sec{The Group of Contactomorphisms}
\definexref{sec-gpofcms}{\the\secno}{section}

\noindent
For a manifold $M$ and a compact subset $K\subset M$, 
let $\dif(M\rel K)$ denote the group of \difmm s 
which become identity when restricted to $K$, i.e. 
$$
\dif(M\rel K)=\{f\in\dif(M) : f\big|_K=\id_K\}.
$$

Let $\xi$ be a \cstr\ on the manifold $M$. 
Denote by $\dif_\xi(M\rel K)$ the group of \cmm s of $(M,\xi)$
coinciding with the identity on $K$:
$$
\dif_\xi(M\rel K)=\{f\in\dif(M\rel K) : f_*\xi=\xi\}.
$$

In \cite{dymara:masters} we have proved that the group $\dif_\xi(S^3\rel K)$ 
for $K$ being a small (closed) ball containing an \ot\ 
disk is not connected; in fact, that $\pi_0$ of this group 
has {\it exactly} two elements. In the present paper we refine 
our argument so as to obtain the complete knowledge of the weak homotopy
type of $\dif_\xi(S^3\rel\Delta)$.


\thm{Theorem}{Let $\xi$ be an \ocstr\ on $S^3$, 
$\Delta$ an \ot\ disk for $\xi$, 
$K\supset\Delta$ a small closed ball containing the \ot\ 
disk. Then the group $\dif_\xi(S^3\rel K)$ 
is weakly homotopy equivalent to $\Omega^4S^2$.}
\definexref{thm-gpofcms}{\here}{thm}


In order to prove this theorem, we first introduce a number of definitions 
and lemmas.

\smallskip
Let $\tdif_\xi(M\rel K)$ be the space of paths in $\dif(M\rel K)$ 
beginning at identity and ending at a \cmm, i.e. 
$$
\eqalign{\tdif_\xi(M\rel K)=\big\{\{f_t, t\in[0,1]\} :
& f_t\in\dif(M\rel K) \hbox{ for all }t, \cr
& f_0=\id, f_1\in\dif_\xi(M\rel K)\big\}.}
$$
Define the map $\tau:\tdif_\xi(M\rel K)\to\dif_\xi(M\rel K)$
as ``taking the endpoint of a path", i.e.~$\tau(\{f_t\})=f_1$. 


\thm{Lemma}{For $K$ a three-ball in $S^3$ the map 
$\tau:\tdif_\xi(S^3\rel K)\to\dif_\xi(S^3\rel K)$ is a weak 
homotopy equivalence.}
\definexref{lem1}{\here}{lem}

\proof
It suffices to check that the fiber $\tau^{-1}(f)$ is contractible 
for any point $f$ in the space $\dif_\xi(S^3\rel K)$.
The whole $\dif(S^3\rel K)$ is contractible 
(Smale conjecture, \cite{hatcher:smale}).
Since $\tau^-1(f)$ consists of all paths  
in $\dif(S^3\rel K)$ with fixed endpoints (joining the identity with $f$), 
it is contractible as well. 
\qed

\bigskip
For a \cstr\ $\xi_K$ on a compact subset $K$ of a manifold $M$,
denote by $\cont(M\rel K,\xi_K)$ the space of all \cstr s on $M$ which 
coincide with $\xi_K$ on $K$. 
Let $\xi_*$ be the map 
$\dif(M\rel K)\to(\cont(M\rel K),\xi\big|_K)$
defined as $\xi_*(f)=f_*(\xi)$.

Note that for $\{f_t\}\in\tdif_\xi(M\rel K)$ the family 
$\{\xi_*(f_t)\}$ is actually an element 
of the space $\Omega(\cont(M\rel K),\xi\big|_K)$,
because both $f_0$ and $f_1$ are \cmm s. Let us call thus induced map
$\Omega\xi_*:\tdif_\xi(M\rel K)\to\Omega\cont(M\rel K)$.

The following lemma is a corollary of Gray's theorem (\ref{thm-gray2}).

\thm{Lemma}{For any  manifold $M$ with a \cstr\ $\xi$, 
and any compact subset $K$,
the map $\Omega\xi_*:\tdif_\xi(M\rel K)\to\Omega\cont(M\rel K,\xi\big|_K)$ 
is a weak homotopy equivalence.}
\definexref{lem2}{\here}{lem}

\medskip
\noindent
{\it Proof:}
We wish to show that the maps induced by $\Omega\xi_*$ on the homotopy groups
are isomorphisms.

Consider a $k$-sphere mapped into $\Omega\cont$, i.e.~a 
family of \cstr s 
$\{\xi^z_t \mid z\in S^k,t\in[0,1]\}$, such that
\item{$\bullet$} $\xi_0^z=\xi_1^z=\id$ for all $z\in S^k$;
\item{$\bullet$} for $\N$ the north pole of $S^k$, 
$\xi_t^\N=\id$ for all $t$ (this condition, as well as any 
analogous condition formulated below, means simply taking the sphere 
with a fixed base point, in compliance with the very definition
of homotopy groups).

This family is parameterized by $S^k{\times}[0,1]$, which is not contractible; 
therefore Gray's theorem does not apply.
Extend the family of \cstr s to the $k+1$-dimensional disk
$D=S^k{\times}[0,1]\cup D^{k+1}{\times}\{0\}$, setting $\xi^z_0=\xi$ for all
$z\in D^{k+1}$. Now the subset of $D$ where we require that the contact structure
coincides with $\xi$ is the heavily shaded area of \ref{fig-param} (including the thick lines).

\midinsert
\centerline{\epsffile{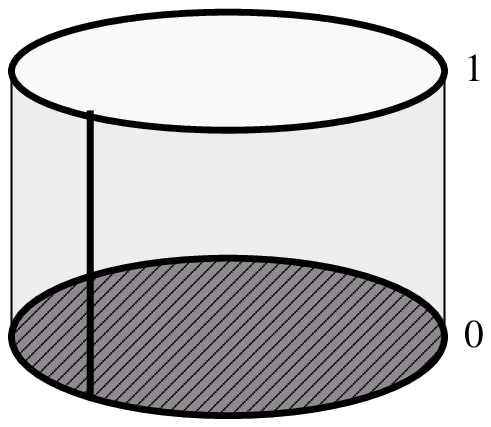}}
\hskip4cm \raise 4mm \hbox{$\cal N$}
\figure{The parameter space for the extended family of \cstr s.}
\definexref{fig-param}{\the\secno.\the\figno}{figure}
\endinsert

By Gray's theorem (\refn{thm-gray2}) the extended family has a trivialization, 
i.e.~a family of projections $f^z_t$ on a fixed fiber.
The relative version of Gray's theorem (\refn{thm-gray2}) allows to assume 
that the desired family of projections extends a certain given family 
of diffeomorphisms 
parameterized by a contractible subset of $D$; we choose the contractible
subset consisting of the disk $D^{k+1}{\times}\{0\}$ and the line segment
$\N{\times}[0,1]$, where the projections can be assumed to be the identity.
Thus we get the family $\{f^z_t\mid(z,t)\in D\}$ such that 
$(f^z_t)_*\xi=\xi^z_t$ for all $z\in S^k$ and all $t$ and 
at the same time satisfying the following conditions:
\item{$\bullet$} $f^z_0=\id$ for all $z\in D^{k+1}$;
\item{$\bullet$} $f^\N_t=\id$ for all $t\in [0,1]$.

Since $\xi^z_1=\xi$ for all $z\in S^k$, the diffeomorphisms 
$f^z_1$ are contactomorphisms.
Therefore for any $z\in S^k$ the family $F^z=\{f^z_t\mid t\in [0,1]\}$
is an element of the group $\tdif_\xi(M\rel K)$, so that any sphere in
$\Omega\cont(M\rel K, \xi\big|_K)$ is indeed the image (under $\Omega\xi_*$) 
of a sphere in $\tdif(M\rel K, \xi)$.
This proves the surjectivity of the homomorphisms of homotopy groups.

\smallskip
In order to prove their injectivity it is enough to show that if the image 
of a sphere is homotopically trivial in $\Omega\cont$ (i.e.~the family
of loops of contactomorphisms parameterized by a sphere can be 
extended to the disk bounded by the sphere) then the same is true for the 
sphere in $\tdif$ itself. 
It turns out that a construction analogous to that employed 
in the first part of the proof allows to pull back to $\tdif$
the disk spanned by the sphere.
Indeed, consider a family of diffeomorphisms 
${\cal F}=\{f^z_t\in \tdif_\xi(M\rel K)\mid z\in S^k, t\in[0,1]\}$,
such that 
\item{$\bullet$} $f^z_0=\id$ for all $z\in S^k$;
\item{$\bullet$} $f^z_1$ is a contactomorphism for all $z\in S^k$;
\item{$\bullet$} for ${\cal N}$ the north pole of $S^k$,
$f^{\cal N}_t=\id$ for all $t$.

\noindent
Then $\Omega\xi_*({\cal F})=\{\xi^z_t\mid z\in S^k, t\in[0,1]\}$
is a sphere in $\Omega\cont$, as described in the first part of 
the proof. In particular,
\item{$\bullet$} $\xi^z_0 = \xi$ for all $z\in S^k$ (because $f^z_0=\id$);
\item{$\bullet$} $\xi^z_1 = \xi$ for all $z\in S^k$ (because
$f^z_1$ is a contactomorphism);
\item{$\bullet$} for ${\cal N}$ the north pole of $S^k$,
$\xi^{\cal N}_t=\xi$ for all $t$.

\noindent
The assumption that $\Omega\xi_*({\cal F})$ is homotopically
trivial in $\Omega\cont$ means that it extends to a family
of contact structures $\{\xi^z_t\mid z\in D^{k+1}, t\in[0,1]\}$
such that $\xi^z_0=\xi^z_1=\xi$ for all $z\in D^{k+1}$.
The subset $D=S^k\times[0,1]\cup D^{k+1}\times\{0\}
\subset D^{k+1}\times[0,1]$ is contractible,
thus again we can apply Gray's theorem (\refn{thm-gray2})
and obtain a family of diffeomorphisms
$\overline{\cal F}=\{{\bar f}^z_t\in\tdif_\xi(M\rel K)\mid
z\in D^{k+1}, t\in[0,1]\}$ such that
\item{$\bullet$} ${\bar f}^z_0=f^z_0=\id$ for all $z\in D^{k+1}$;
\item{$\bullet$} ${\bar f}^z_t=f^z_t$ for all $(z,t)\in D$.

\noindent
Therefore $\overline{\cal F}$ gives a $(k+1)$-disk in 
$\tdif_\xi(M\rel K)$ spanned by the original sphere ${\cal F}$.

This concludes the proof.
\qed

\medskip

Now, fix a plane field $\zeta_K$ on a compact subset $K\subset M$. 
Let ${\it Distr}(M\rel K,\zeta_K)$ denote the space of all distributions 
on $M$ coinciding with $\zeta_K$ on $M$.

The following two lemmas are easy observations.

\thm{Lemma}{For any co-oriented, parallelized distribution $\zeta$ there is a 
homeomorphism 
${\it Distr}(M\rel K,\zeta\big|_K)\to{\it Map}\bigl((M,K)\to(S^2,{\cal N}):
\xi\mapsto G_\xi\bigr)$, where the map $G_\xi$ is defined using the 
parallelization of $M$ induced by that of $\zeta$.}
\definexref{lem3}{\here}{lem}

\thm{Lemma}{For $K$ a three-ball in $S^3$ the space
${\it Map}\bigl((S^3,K)\to(S^2,{\cal N})\bigr)$ is homeomorphic to
$\Omega^3S^2$.}
\definexref{lem4}{\here}{lem}

We conclude this section with the proof of theorem 
formulated at its beginning.
\bigskip
\noindent
{\it Proof of \ref{thm-gpofcms}:}

The subset $K=\overline{U(\Delta)}$ is a three-ball; therefore by \ref{lem1} 
$\dif_\xi(S^3\rel K)$ is \heq t to 
$\tdif_\xi(S^3\rel K)$.

By \ref{lem2} with $M=S^3$, the space $\tdif_\xi(S^3\rel K)$ is \heq t to 
$\Omega(\cont(S^3\rel K),\xi\big|_K)$.

By Eliashberg's theorem with $M=S^3$, the space
$\cont(S^3\rel K,\xi\big|_K)$ is \heq t to 
${\it Distr}(S^3\rel K,\xi\big|_K)$. 
The \heq ce induces a \heq ce on the loop spaces, therefore also 
$\Omega\cont(S^3\rel K,\xi\big|_K)$ 
is \heq t to $\Omega{\it Distr}(S^3\rel K,\xi\big|_K)$. 

By \ref{lem3} with $M=S^3$, ${\it Distr}(S^3\rel K,\zeta\big|_K)$ 
is homeomorphic (therefore \heq t) to 
${\it Map}\bigl((S^3,K)\to(S^2,{\cal N})\bigr)$. Thus also 
$\Omega{\it Distr}(S^3\rel K,\zeta\big|_K)$ is homeomorphic 
to $\Omega{\it Map}\bigl((S^3,K)\to(S^2,{\cal N})\bigr)$.

By \ref{lem4}, $\Omega{\it Map}\bigl((S^3,K)\to(S^2,{\cal N})\bigr)$ 
is homeomorphic to $\Omega(\Omega^3S^2)=\Omega^4S^2$.
\qed

\sec{Discussion}

In the previous section we studied 
the group $\dif_\xi(S^3\rel K)$ of \cmm s fixing (pointwise)
a \nbd\ of an \ot\ disk, rather than the full group of \cmm s
$\dif_\xi(S^3)$. Such is the price we had to pay 
for the knowledge of the full
homotopy type of the group of \cmm s (due to the restrictions 
placed by the assumptions of 
\ref{thm-eliashberg}).
The only known result about the full group $\dif_\xi(S^3)$
(for $\xi$ \ot) is Theorem~2.2.1 in \cite{dymara:masters}: 
it is not connected.
\medskip

Fix an \ot\ disk $\Delta$ and its small open \nbd $K$.
By a {\it contact embedding} of $K$ into 
the contact manifold $(S^3,\xi)$
we will understand a map $f:K\hookrightarrow S^3$
such that $f:K\to f(K)$ is a \cmm.
Two such embeddings are {\it contact isotopic}\/
(or simply {\it isotopic}), if they are homotopic through
contact embeddings.
We will say that two \ot\ disks are (contact) isotopic if they are 
images of $\Delta$ by two isotopic embeddings of $K$.
Denote the space of contact embeddings of 
the \nbd\ $K$ of an \ot\ disk into $S^3$
by ${\it Emb}_\xi(K{\hookrightarrow}S^3)$
and consider the fibration
$$
\dif_\xi(S^3\rel K) \to \dif_\xi(S^3) \to 
{\it Emb}_\xi(K{\hookrightarrow}S^3).
$$
One can see that what we are lacking to grasp
the structure of the full group of \cmm s
is exactly understanding the space of contact embeddings of
the \ot\ disk;
\ref{prop-disjot} gives a partial result in this 
direction.
Two lemmas are needed in the proof;
the first one essentially has been contained in \dy,
but here for the first time we state it explicitly.

\thm{Lemma}{Let $K$, $L$ be compact subsets of a closed 
contact manifold $(M,\xi)$ such that $M\setminus K$ is connected,
$\Delta$ an \ot\ disk in $M\setminus (K\cup L)$.
Let $f:M\to M$ be a diffeomorphism which fixes a \nbd\ of $\Delta$,
moves $K$ onto $L$, and becomes a \cmm\ when restricted to a \nbd\ of 
$K$. Assume that the \cstr s $\xi$ and $f_*\xi$ are 
homotopic as plane fields on $M\rel(K\cup\Delta)$. Then
there is a \cmm\ $F:M\to M$ moving $K$ onto $L$.}
\definexref{lem-cmm}{\here}{lemma}

\proof
By Eliashberg's theorem (\refn{thm-eliashberg}) $\xi$ and $f_*\xi$ are
isotopic. By Gray's theorem (\refn{thm-gray2}) they are also 
isomorphic, i.e.~there is a diffeomorphism $g$ such that 
$g_*\xi=f_*\xi$, $g\bigm|_K=\id$. Thus $F=f\circ g^{\-1}$ is a \cmm\
and $F(K)=L$ as desired.
\qed

The assumptions of \ref{lem-cmm} can be rephrased in the language of 
Gauss maps with respect to the parallelization induced by the 
parallelization of $\xi$ as follows:
there is a homotopy $G_t:M\to S^2$ such that $G_0=G_\xi$ is constant,
$G_t\bigm|_{K\cup\Delta}$ is constant for all $t$ 
and $G_1=G_{f_*\xi}$.

\thm{Corollary}{Let $K\subset \r^3$ be a contractible compact set,
$\xi_K$ a \cstr\ on a \nbd\ $U_K$ of $K$. Furthermore, 
let $j_1,j_2$ be two embeddings of $U_K$ into a contact manifold $M$
which are \cmm s onto their respective images. 
If there is an \ot\ disk $\Delta$ in 
$M\setminus\bigl(j_1(U_K)\cup j_2(U_K)\bigr)$, then
there exists a \cmm\ of $M$ moving $j_1(K)$ onto $j_2(K)$.}
\definexref{coro-embs}{\here}{coro}

\proof
The map $j_2^{-1}\circ j_1$ can be extended to a diffeomorphism 
$J:M\to M$ coinciding with the identity outside an open ball
$B\supset j_1(K)\cup j_2(K)$, $B\cap\Delta=\emptyset$.
Since the group $\dif(B\rel\partial B)$ is connected, $J$
must be isotopic to the identity; denote this isotopy
by $J_t:M\to M$, so that we have $J_t\bigm|_{M\setminus B}=\id,
J_0=\id,  J_1=J$. Now, consider the homotopy of plane fields 
given by $(J_t)_*\xi$, and the associated homotopy of Gauss maps
$G_t:M\to S^2$. We see that $G_0$, as well as $G_t\bigm|_{M\setminus B}$
and $G_1\bigm|_{j_1(K)}$ , are constant.
Since $K$ is contractible, the map ${\cal G}:K\times[0,1]\to S^2:
(x,t)\mapsto G_t\circ J_t(x)$ is homotopic (rel $K\times\{0,1\}$)
to the constant map. This homotopy can be extended to the whole $M$,
yielding (by the $\xi\leftrightarrow G_\xi$ identification)
the plane field homotopy satisfying the assumptions of \ref{lem-cmm}.
\qed

\thm{Lemma}{Let $L\subset S^3$ be a compact set,
$f:S^3\to S^3$ a contactomorphism coinciding wit the identity
on (a \nbd\ of) an \ot\ disk $\Delta$ disjoint with $L$.
Then $L$ and $f(L)$ are ambientally isotopic, 
i.e.~there exists a family of \cmm s 
$\{f_t : t\in[0,1]\}, f_0=\id, f_1(L)=f(L)$.}
\definexref{lem-ambient}{\here}{lemma}

\proof
As before, let $K=\overline{U(\Delta)}$. The weak homotopy
equivalence $\dif_\xi(S^3\rel K)\to \Omega^4S^2$ postulated by 
\ref{thm-gpofcms} induces a map 
${\cal H}:\dif_\xi(S^3\rel K)\to\pi^4S^2\simeq\z^2$.
Choose a \cmm\ $g$ coinciding with the identity 
outside a small ball disjoint with $K$ and $L$ and 
such that ${\cal H}(g)={\cal H}(f)^{-1}$. Then
$g\circ f$, since it lies in the kernel of ${\cal H}$,
is isotopic to the identity through \cmm s.
This isotopy provides the desired family of \cmm s.
\qed

Note that ambientally isotopic \ot\ disks are necessarily isotopic.

\thm{Proposition}{Two disjoint \ot\ disks are isotopic.}
\definexref{prop-disjot}{\here}{propo}

\proof
Let $\Delta_1$ and $\Delta_2$ be the disjoint \ot\ disks.
A parallel copy of $\Delta_1$ (close enough to it) is an \ot\ disk
disjoint with both $\Delta_1$ and $\Delta_2$.
By \ref{coro-embs} there is a contactomorphism moving $\Delta_1$ onto $\Delta_2$.
Then \ref{lem-ambient} implies that $\Delta_1$ is isotopic to $\Delta_2$.
\qed
\bigskip

No actual example of two non-isotopic \ot\ disks is known to the author.
One may, however, try contemplating the following construction.

Let $T$ be the solid torus 
$$
\{(r,\theta,z) : r<\pi+\eps\}/_{(r,\theta,z)\sim(r,\theta,z+2\pi)}
$$
with the contact structure
$$
\xi=\ker\bigl(\cos r - \varrho(r)\,\sin r\,d\theta\bigr),
$$
as in \ref{sec-prelim}.
Identify two copies of $T$ along the subset 
$$
N=\{(r,\theta,z) : \hbox{\frac$\pi$/2}-\eps<r<\pi+\eps\}
$$
via the diffeomorphism 
$$
\phi:(r,\theta,z)\mapsto(\hbox{\frac3$\pi$/2} -r,z,\theta).
$$
Topologically, the manifold we obtain is the three-sphere; one should think
about the genus~1 Heegaard splitting of $S^3$ 
with both solid tori slightly enlarged, so that the identification,
rather than along the boundary, is performed along its tubular \nbd.
We shall denote the two inclusions of $T$ into $S^3$ by $\psi_1$ and $\psi_2$;
on $N$ we have $\phi=\psi_2^{-1}\circ\psi_1$.

Moreover, if we choose the function $\varrho$ so that $\varrho(r)=1$ for $r\in(\frac$\pi$/2-\eps,\pi+\eps)$, then 
$\phi$ is a \cmm, so the resulting sphere carries a \cstr\ pushed forward 
from the solid tori.
This \cstr\ is overtwisted, with two obvious families of \ot\ disks of the form
$\psi_i\bigl(\{(r,\theta,z) : r\leq \pi, z=z_0\}\bigr)$  for $i=1,2$.
Note that the formula 
$\Psi(x)=
\cases{\psi_2\circ\psi_1^-1(x) & for $x\in\psi_1(T)$, \cr 
\psi_1\circ\psi_2^-1(x) & for $x\in\psi_2(T)$}$
defines a \cmm\ of $S^3$ exchanging these two families.

Consider two disks, one of each of the families, say
$\Delta_1=\psi_1\bigl(\{(r,\theta,z) : r\leq \pi, z=0\}\bigr)$ and
$\Delta_2=\psi_2\bigl(\{(r,\theta,z) : r\leq \pi, z=0\}\bigr)$.
It is easy to see that $\Delta_1$ and $\Delta_2$ intersect along an interval,
and not much more difficult---that $\xi$ restricted to 
$S^3\setminus(\Delta_1\cup\Delta_2)$ is tight\numberedfootnote{For
a detailed proof that 
even $S^3\setminus(\partial\Delta_1\cup\partial\Delta_2)$ is tight 
see \cite{dymara:knots}.}; 
therefore assumptions of neither \ref{coro-embs} nor \ref{prop-disjot} are satisfied.

\thm{Question}{Are $\Delta_1$ and $\Delta_2$ isotopic?}

\sec{References}

\bibliographystyle{siam}
\bibliography{contact,lanl}

\bye

\medskip
The tight-vs-overtwisted dichotomy has unexpectedly deep 
consequences. Methods useful in
studying \ot\ contact manifolds 
have very little in common with those appropriate
for tight ones. The central result in \ot\ contact topology is
the following

\sec{Overtwisted Contact Structures}
\definexref{sec-ocstr}{\the\secno}{section}

Each of the \cstr s described in 
\ref{ex-zetan} belongs to one of exactly three different 
isomorphism classes. It is easy to notice that, 
considered as plane fields, all $\zeta_n$'s
with $n$ odd realize one homotopy type, those with $n$ even---another one.
But in addition to that, $\zeta_0$ can be proved to be different
from $\zeta_2$ ($\simeq\zeta_4\simeq\ldots$) using the notion of 
\ot ness (introduced in \cite{bennequin}).

\medskip
A way to describe explicitly each of the infinitely many \ocstr s
on $S^3$ is provided by Lutz diagrams (see \cite{lutz:twists}; cf.~also 
\cite{dymara:masters}). Lutz considers \cstr s invariant under 
the Hopf action of $S^1$ on $S^3$ (each homotopy type can be 
realized in this manner) and in general position.
The set $\Sigma$ of fibers tangent to the
\cstr\ forms a codimension 1 submanifold in the 
quotient $S^2$, i.e.~a family of circles.
Then we mark the components of $S^2\setminus\Sigma$ with $+$ or $-$,
depending whether the fiber intersects the contact planes 
in the direction of their co-orientation or in the opposite direction.
These markings on $S^2$---the family of circles and the signs---form
the Lutz diagram of the given \cstr.

In particular, the Lutz diagram of the structure $\zeta_0$ 
of \ref{ex-zetan} is empty: there are no fibres tangent
to the contact planes. The diagram of the \cstr{} $\zeta_1$
consists of a single circle, dividing
$S^2$ into two regions, one of which is marked $+$ and the other $-$.

It is easy to find, using \ref{cube}, \cstr s with Lutz diagrams 
consisting of any number of concentric circles.

\rmk{Example}
The \cstr\ $\zeta_n$ on $S^3$ defined in \ref{ex-zetan}
has the Lutz diagram made of $n$ concentric circles.

\midinsert
\figure{The Lutz diagram of $\zeta_n$.}
\endinsert

\smallskip
The homotopy type $\eta(\xi)$ of a \cstr\ $\xi$ can be read from 
its Lutz diagram $\Sigma(\xi)$ as follows.
Consider the circles as lying in $\r^2$ rather than $S^2$
by removing a point from one of the regions marked $+$
\numberedfootnote{This is
the only artificial choice made in this construction.
It corresponds to the fact that 
the actual identification of $\pi_0\cont(S^3)$ with $\pi_3S^2$ 
depends on the choice of the parallelization.}.
To each circle $c$ assign a number 
$$\eta(c)=\cases{1 & if $c$ has $+$ inside and $-$ outside, \cr
                -1 & otherwise.}$$
Then $\eta(\xi)=\sum_{c\in\Sigma(\xi)}\eta(c)\in{\bf Z}\simeq\pi_3S^2$.

\rmk{Example}
For $\zeta_n$ of \ref{ex-zetan} we have $\eta(\zeta_n)=\cases{
-1 & for $n$ odd,  \cr 
0 & for $n$ even.}$

\smallskip
One can find \cstr s of other homotopy types using so called Lutz twists.
This operation is a (topologically trivial) surgery on a closed curve
transversal to the \cstr. Replace the \cstr\ on its tubular \nbd\ 
with one contactomorphic to the \cstr\ $\zeta_1$
restricted to $\pr_{\cal C}([0,\frac2/3+\eps]\times[0,1]\times[0,1])$,
i.e.~a \nbd\ of the sum of a family of \ot\ disks indexed by a circle
(such a family is called an {\it \ot\ tube}).
Thus, the \cstr\ obtained by a Lutz twist is always \ot.

\rmk{Proposition}
Every \ocstr\ arises from a Lutz twist.
\definexref{prop-fromlutz}{\here}{propo}

\proof
Notice that 
Lutz twist done twice along the same curve leaves
the plane field homotopy type of the \cstr\ unchanged.
Denote by $LT_k(\xi)$ the \cstr\ obtained from $\xi$ 
by the Lutz twist along $k$. If $\xi$ is \ot,
we have $\xi=LT_k(LT_k(\xi))$ for any transversal knot $k$.
\qed

Caution: this does not mean that every \ocstr\ arises from a Lutz 
twist {\it applied to a tight \cstr}. It is the case 
if and only if there is a tight \cstr\ in the appropriate 
plane field homotopy class.
On $S^3$, for example, there is a unique tight \cstr, namely
$\zeta_0$ of \ref{ex-zetan};
and therefore $\zeta_1$ is the only \cstr\ which arises from
the Lutz twist performed on a tight \cstr.
However, starting with $\zeta_0$ and performing consecutive
Lutz twists one can eventually reach every \cstr 
(and generally, if a manifold possesses a tight \cstr,
then any \ocstr\ can be reached from it by a number of Lutz twists;
\cite{lutz:twists}). 

One can formulate a sort of relative version of
\ref{prop-fromlutz}, providing answers to 
Question~35 and (partially) Question~36 of \cite{problems}.

\rmk{Proposition}
An \ot\ disk in can always be completed to 
an \ot\ tube produced by a Lutz twist.

\proof
Let $\Delta$ be an \ot\ disk in $(M,\xi)$.
Perform a double Lutz twist along a curve in $M\setminus\Delta$.
The obtained contact structure $\xi''$ is isomorphic to the original one,
and the isomorphism can be chosen so as to coincide with 
the identity on $\Delta$. In $(M,\xi'')$ the disk $\Delta$ 
can be proven (using Eliashberg's theorem and the existence of 
an \ot\ disk disjoint with both $\Delta$ and the support 
of the Lutz twist)
to be isotopic to an \ot\ disk in the newly formed 
\ot\ tube obtained from the Lutz twist.
\qed

let us employ the following way of looking at $S^3$.

\midinsert
\figure{Gluing of the cube (\ref{cube}).}
\definexref{fig-cube}{\the\secno.\the\figno}{figure}
\bigskip
\endinsert

\rmk{Construction} \definexref{cube}{\here}{construction}
Consider the cube  ${\cal C}=[0,1]\times[0,1]\times[0,1]$. 
Identify two points $(x,y,z)$ and $(x',y',z')$ in $\cal C$
if any of the following condition holds:\numberedfootnote{In 
more informal language we would say: glue together 
the left and right faces of the cube, 
and also the top and bottom faces; then contract to points all 
horizontal lines on the front face, and also all vertical lines 
on the back face. Cf.~\ref{fig-cube}.}
\item{(a)} $x=x', y=y', z=0, z'=1$, or
\item{(b)} $x=x', y=0, y'=1, z=z'$, or
\item{(c)} $x=x'=0, z=z'$, or
\item{(d)} $x=x'=1, y=y'$.

Then quotient is homeomorphic to $S^3$.
We will denote the projection map ${\cal C}\to S^3$
by $\pr_{\cal C}$.

\medskip

This construction provides a system of coordinates on the 3-sphere, 
which can be used to write explicitly 1-forms defining certain \cstr s.

\rmk{Example} \definexref{ex-zetan}{\here}{exple}
For any integer $n\geq 0$, let the \cstr\ $\zeta_n$ on $S^3$ 
be the kernel of the 1-form 
$\alpha_n=\cos\bigl( (n+\frac1/2)\pi x\bigr)\,dz
-\sin\bigl( (n+\frac1/2) \pi x\bigr)\,dy$.
\smallskip


Note that for each $n$ the form $\alpha_n=dz$ on the front face 
$\{x=0\}$ and $\alpha_n=\pm dy$ on the back face $\{x=1\}$.
Thus all $\zeta_n$'s are indeed well-defined \cstr s on $S^3$.
Now let us have a closer look at one of them.

\rmk{Example}
\definexref{ex-stdonsphere}{\here}{exple}
For $n=0$ we get the 
\cstr\ $\zeta_0$ known as the standard \cstr\ on $S^3$. 
It can be described independently of \ref{cube}
as the field of planes perpendicular 
to the fibers of the Hopf bundle. 
Removing a point from $(S^3,\zeta_0)$
yields a \cstr\ on $\r^3$ isomorphic to to the one defined
in \ref{ex-std} (which hopefully justifies the
shameless abuse of notation which we commit
referring to both of them 
by the same symbol $\zeta_0$).

\bye